\documentclass[12pt,twoside]{article}
\usepackage[english]{babel}
\usepackage[latin1]{inputenc}
\usepackage{amsmath}
\usepackage{amssymb,amsfonts}
\usepackage{graphicx}
\usepackage{times,amssymb,amscd}

\numberwithin{equation}{section}

\newcommand{\bA}{\mathbf{A}}

\newcommand{\bE}{\mathbf{E}}

\newcommand{\bH}{\mathbf{H}}

\newcommand{\bL}{\mathbf{L}}

\newcommand{\bR}{\mathbf{R}}
\newcommand{\bS}{\mathbf{S}}
\newcommand{\bV}{\mathbf{V}}

\newcommand{\be}{\mathbf{e}}

\newcommand{\br}{\mathbf{r}}
\newcommand{\bx}{\mathbf{x}}
\newcommand{\by}{\mathbf{y}}

\newcommand{\bT}{\mathbf{T}}

\newcommand{\bt}{\mathbf{t}}

\newcommand{\BV}{\boldsymbol{V}}

\newcommand{\Be}{\boldsymbol{e}}
\newcommand{\Bu}{\boldsymbol{u}}
\newcommand{\Bv}{\boldsymbol{v}}

\newcommand{\cP}{\mathcal{P}}
\newcommand{\cS}{\mathcal{S}}

\newcommand{\cC}{\mathcal{C}}

\newcommand{\EUC}{\mathbf E^3}
\newcommand{\SPH}{\bS^3}
\newcommand{\HYP}{\bH^3}
\newcommand{\SXR}{\bS^2\!\times\!\bR}
\newcommand{\HXR}{\bH^2\!\times\!\bR}
\newcommand{\SLR}{\widetilde{\bS\bL_2\bR}}
\newcommand{\NIL}{\mathbf{Nil}}
\newcommand{\SOL}{\mathbf{Sol}}

\begin{document}
\pagestyle{myheadings}
\markboth{\centerline{Jen\H o Szirmai}}
{Cylinder-like Pappus's hexagon theorem in $\NIL$ geometry}
\title
{Cylinder-like Pappus's hexagon theorem in $\NIL$ geometry
}

\author{Jen\H o Szirmai \\
\normalsize Department of Algebra and Geometry, Institute of Mathematics,\\
\normalsize Budapest University of Technology and Economics, \\
\normalsize M\H uegyetem rkp. 3., H-1111 Budapest, Hungary \\
\normalsize szirmai@math.bme.hu
\date{\normalsize{\today}}}
\maketitle
\begin{abstract}
In this paper we deal with $\NIL$ geometry, 
in whose projective model the geodesic curves are helix-like and fit onto geodesic $\NIL$ cylinders of revolution with fibrum axes. 
In this paper we investigate relations for geodesic cylinders and thus also geodesic 
curves, which lead to an analogous result to Pappus's hexagon theorem and provide important information about the structure of the considered space.
\end{abstract}
\newtheorem{Theorem}{Theorem}[section]
\newtheorem{corollary}[Theorem]{Corollary}
\newtheorem{lemma}[Theorem]{Lemma}
\newtheorem{exmple}[Theorem]{Example}
\newtheorem{definition}[Theorem]{Definition}
\newtheorem{rmrk}[Theorem]{Remark}
\newtheorem{proposition}[Theorem]{Proposition}
\newenvironment{remark}{\begin{rmrk}\normalfont}{\end{rmrk}}
\newenvironment{example}{\begin{exmple}\normalfont}{\end{exmple}}
\newenvironment{acknowledgement}{Acknowledgement}

\section{Introduction}
\label{section1}
The internal structure of non-Euclidean geometries, and thus Thurston geometries, also differs 
significantly from that of the usual geometries with constant curvature, 
but each of them contains various aspects of these.

The analysis of the structure of Thurston geometries 
was greatly facilitated by the description of the 
projective models of Thurston geometries, that we will use in the following (see \cite{M97}). 

In the study of the internal structures of Thurston spaces, it is worth distinguishing 
two significantly different directions related to the translation distance or geodesic distance.

We can introduce in a natural way (see \cite{M97})
translation mappings any point to any other point. Consider a unit tangent
vector at the origin. Translations carry this vector to a tangent vector any
other point.
If a curve $t\rightarrow (x(t),y(t),z(t))$ has just the translated vector as its tangent vector at
each point, then the curve is called a {\it translation curve}. This assumption leads to
a system of first order differential equations. Thus, translation curves are simpler
than geodesics and differ from them in $\NIL$, $\SLR$ and $\SOL$ geometries. Moreover, in the $\NIL$, $\SLR$ and $\SOL$ geometries,
translation curves are in many ways more natural than geodesics.
In $\EUC$, $\SPH$, $\HYP$, $\SXR$ and $\HXR$ geometries, the translation and geodesic curves 
coincide
with each other. Therefore, we
distinguish two different distance functions: $d^g$ is the usual geodesic distance function,
and $d^t$ is the translation distance function. So we obtain
two types of curves, triangles, bisector surfaces (and two types of the corresponding Dirichlet-Voronoi cells) etc. by the two different distance functions,
but {\it{in the present paper we consider only the geodesic case}}. 

We examined several interesting problems related to the Thurston geometries 
with non-constant curvature in papers \cite{Cs-Sz23}, \cite{Cs-Sz25}, \cite{MSzV}, \cite{MSz12}, \cite{Sz07}, 
\cite{Sz14-1}, \cite{Sz14-2}, \cite{Sz18}, \cite{Sz23-2}, \cite{Sz24}. 

Here I would only highlight the results related to the Menelaus' and Ceva's theorems that are particularly significant
in understanding external triangle properties and collinearity conditions.
While its formulation is well established in geometries of constant curvature, its behavior
in other geometries requires careful investigation. In \cite{Sz22-2} and \cite{Sz23-1} we generalized 
them to $\SXR$, $\HXR$ and $\NIL$ spaces. In \cite{Sz26} we continued to investigate this issue, 
generalized and proved the theorems of Menelaus' and Ceva's for {\it translation triangles} in $\NIL$, $\SLR$ and $\SOL$ spaces. 

Continuing these investigations in \cite{Sz26-1} we analyzed Pappus's theorem and Desargues's theorem in Thurston geometries of non-constant curvatures 
in {\it translation case}. However, in similar question related to geodesic curves, the investigation of the geodesic 
Pappus's theorem can only be carried out in certain spaces due to the 
complexity of the geodesic equations. For example, in $\SOL$ geometry, the exact equation of the geodesic curves cannot even be determined.

{\it This paper aims to fill this gap by analyzing cylinder-like Pappus's theorem in $\NIL$ geometry.  
Here we formulate a statement analogous to Pappus's theorem for $\NIL$ cylinders and geodesics, using Euclidean inversion.}

This may also be interesting because its variants can be extended to other Thurston geometries, but in this work we study only $\NIL$ space.
By establishing this theorem provides new insights into the relationships between Euclidean
and non-Euclidean geometries.

\section{Projective model of $\NIL$ geometry}
E. {Moln\'ar} has shown in \cite{M97}, that the homogeneous 3-spaces
have a unified interpretation in the projective 3-sphere $\mathcal{PS}^3(\bV^4,\BV_4, \mathbf{R})$. 
In our work we shall use this projective model of $\NIL$ geometry. 
The Cartesian homogeneous coordinate simplex is given by $E_0(\be_0)$,$E_1^{\infty}(\be_1)$,$E_2^{\infty}(\be_2)$,
$E_3^{\infty}(\be_3)$, $(\{\be_i\}\subset \bV^4$ \ and $\text{with the unit point}$ $E(\be = \be_0 + \be_1 + \be_2 + \be_3 ))$. 
Moreover, $\by=c\bx$ with $0<c\in \mathbf{R}$ (or $c\in\mathbf{R}\setminus\{0\})$
defines a point $(\bx)=(\by)$ of the projective 3-sphere $\cP \cS^3$ (or that of the projective space $\cP^3$ where opposite rays
$(\bx)$ and $(-\bx)$ are identified). 
The dual system $\{(\Be^i)\}\subset \BV_4$ describes the simplex planes, especially the plane at infinity 
$(\Be^0)=E_1^{\infty}E_2^{\infty}E_3^{\infty}$, and generally, $\Bv=\Bu\frac{1}{c}$ defines a plane $(\Bu)=(\Bv)$ of $\cP \cS^3$
(or that of $\cP^3$). Thus $0=\bx\Bu=\by\Bv$ defines the incidence of point $(\bx)=(\by)$ and plane
$(\Bu)=(\Bv)$, as $(\bx) \text{I} (\Bu)$ also denote it. Thus, {$\NIL$} can be visualized in the affine 3-space $\bA^3$
(so in $\bE^3$) as well.
\subsection{Geodesic curves and spheres in $\NIL$ space}
In this section we recall the important notions and results from the papers \cite{M97}, \cite{PSSz10}, \cite{Sz07}, \cite{Sz23-1}.

$\NIL$ geometry is a homogeneous 3-space derived from the famous real matrix group $\mathbf{L(R)}$, used by W.~Heisenberg in his electro-magnetic studies.
The Lie theory with the method of projective geometry makes possible to describe this topic.

The left (row-column) multiplication of Heisenberg matrices
     \begin{equation}
     \begin{gathered}
     \begin{pmatrix}
         1&x&z \\
         0&1&y \\
         0&0&1 \\
       \end{pmatrix}
       \begin{pmatrix}
         1&a&c \\
         0&1&b \\
         0&0&1 \\
       \end{pmatrix}
       =\begin{pmatrix}
         1&a+x&c+xb+z \\
         0&1&b+y \\
         0&0&1 \\
       \end{pmatrix}
      \end{gathered} \tag{2.1}
     \end{equation}
defines "translations" $\mathbf{L}(\mathbf{R})= \{(x,y,z): x,~y,~z\in \mathbf{R} \}$ on the points of
$\NIL= \{(a,b,c):a,~b,~c \in \mathbf{R}\}$.
These translations are not commutative, in general. The matrices $\mathbf{K}(z) \vartriangleleft \mathbf{L}$ of the form
     \begin{equation}
     \begin{gathered}
       \mathbf{K}(z) \ni
       \begin{pmatrix}
         1&0&z \\
         0&1&0 \\
         0&0&1 \\
       \end{pmatrix}
       \mapsto (0,0,z)
      \end{gathered}\tag{2.2}
     \end{equation}
constitute the one parametric centre, i.e., each of its elements commutes with all elements of $\mathbf{L}$.
The elements of $\mathbf{K}$ are called {\it fibre translations}. $\NIL$ geometry of the Heisenberg group can be projectively
(affinely) interpreted by the "right translations"
on points as the matrix formula
     \begin{equation}
     \begin{gathered}
       (1;a,b,c) \to (1;a,b,c)
       \begin{pmatrix}
         1&x&y&z \\
         0&1&0&0 \\
         0&0&1&x \\
         0&0&0&1 \\
       \end{pmatrix}
       =(1;x+a,y+b,z+bx+c)
      \end{gathered} \tag{2.3}
     \end{equation}
shows, according to (2.1). Here we consider $\mathbf{L}$ as projective collineation
group with right actions in homogeneous coordinates.

In this context E. Moln\'ar \cite{M97} has derived the well-known infinitesimal arc-length-square, invariant under translations $\bL$ at any point of $\NIL$ as follows
\begin{equation}
   \begin{gathered}
      (dx)^2+(dy)^2+(-xdy+dz)^2=\\
      (dx)^2+(1+x^2)(dy)^2-2x(dy)(dz)+(dz)^2=:(ds)^2
       \end{gathered} \tag{2.4}
     \end{equation}
Hence we get the symmetric metric tensor field $g$ on $\NIL$ by components $g_{ij}$, furthermore, its inverse:
\begin{equation}
   \begin{gathered}
       g_{ij}:=
       \begin{pmatrix}
         1&0&0 \\
         0&1+x^2&-x \\
         0&-x&1 \\
         \end{pmatrix},  \quad  g^{ij}:=
       \begin{pmatrix}
         1&0&0 \\
         0&1&x \\
         0&x&1+x^2 \\
         \end{pmatrix} \\
         \text{with} \ \det(g_{ij})=1.
        \end{gathered} \tag{2.5}
     \end{equation}
The translation group $\mathbf{L}$ defined by formula (2.3) can be extended to a larger group $\mathbf{G}$ of collineations,
preserving the fibering, that will be equivalent to the (orientation preserving) isometry group of $\NIL$.
In \cite{M06} E.~Moln\'ar has shown that
a rotation trough angle $\omega$
about the $z$-axis at the origin, as isometry of $\NIL$, keeping invariant the Riemann
metric everywhere, will be a quadratic mapping in $x,y$ to $z$-image $\overline{z}$ as follows:
     \begin{equation}
     \begin{gathered}
       \br(O,\omega):(1;x,y,z) \to (1;\overline{x},\overline{y},\overline{z}); \\
       \overline{x}=x\cos{\omega}-y\sin{\omega}, \ \ \overline{y}=x\sin{\omega}+y\cos{\omega}, \\
       \overline{z}=z-\frac{1}{2}xy+\frac{1}{4}(x^2-y^2)\sin{2\omega}+\frac{1}{2}xy\cos{2\omega}.
      \end{gathered} \tag{2.6}
     \end{equation}
This rotation formula, however, is conjugate by the quadratic mapping 
     \begin{equation}
     \begin{gathered}
       \mathcal{M}:~x \to x'=x, \ \ y \to y'=y, \ \ z \to z'=z-\frac{1}{2}xy  \ \ \text{to} \\
       (1;x',y',z') \to (1;x',y',z')
       \begin{pmatrix}
         1&0&0&0 \\
         0&\cos{\omega}&\sin{\omega}&0 \\
         0&-\sin{\omega}&\cos{\omega}&0 \\
         0&0&0&1 \\
       \end{pmatrix}
       =(1;x",y",z"), \\
       \text{with} \ \ x" \to \overline{x}=x", \ \ y" \to \overline{y}=y", \ \ z" \to \overline{z}=z"+\frac{1}{2}x"y",
      \end{gathered} \tag{2.7}
     \end{equation}
i.e. to the linear rotation formula. This quadratic conjugacy modifies the $\NIL$ translations in (2.3), as well.
This can also be characterized by the following important classification theorem.
\begin{Theorem}[E.~Moln\'ar \cite{M06}]
\begin{enumerate}
\item Any group of $\NIL$ isometries, containing a 3-dimensional translation lattice,
is conjugate by the quadratic mapping in (2.5) to an affine group of the affine (or Euclidean) space $\bA^3=\EUC$
whose projection onto the (x,y) plane is an isometry group of $\bE^2$. Such an affine group preserves a plane
$\to$ point null-polarity.
\item Of course, the involutive line reflection about the $y$ axis
     \begin{equation}
     \begin{gathered}
       (1;x,y,z) \to (1;-x,y,-z),
      \end{gathered} \notag
     \end{equation}
preserving the Riemann metric, and its conjugates by the above isometries in {$1$} (those of the identity component)
are also {$\NIL$}-isometries. There does not exist orientation reversing $\NIL$-isometry.
\end{enumerate}
\end{Theorem}
\begin{rmrk}
We obtain a new projective model of $\NIL$ geometry from the above projective model, derived by the above quadratic mapping $\mathcal{M}$.
This is the {\it linearized model of $\NIL$ space} (see \cite{M06}) that seems to be more advantageous to the future investigations. 
But we remain in the classical so called Heisenberg model in this paper.
\end{rmrk}
\subsection{Geodesic curves and spheres} \label{subsection2}
The geodesic curves of the $\NIL$ geometry are generally defined as having locally minimal arc length between their any two (near enough) points.
The equation systems of the parametrized geodesic curves $g^\NIL(x(t),y(t),z(t))$ in our model (now by (2.4)) can be determined by the
Levy-Civita theory of Riemann geometry.
We can assume, that the starting point of a geodesic curve is the origin because we can transform a curve into an
arbitrary starting point by translation (2.1);
\begin{equation}
\begin{gathered}
        x(0)=y(0)=z(0)=0; \ \ \dot{x}(0)=c \cos{\alpha}, \ \dot{y}(0)=c \sin{\alpha}, \\ \dot{z}(0)=w; \ - \pi \leq \alpha \leq \pi. \notag
\end{gathered}
\end{equation}
The arc length parameter $s$ is introduced by
\begin{equation}
 s=\sqrt{c^2+w^2} \cdot t, \ \text{where} \ w=\sin{\theta}, \ c=\cos{\theta}, \ -\frac{\pi}{2}\le \theta \le \frac{\pi}{2}, \notag
\end{equation}
i.e. unit velocity can be assumed.

The equation systems of a helix-like geodesic curves $g^\NIL(x(t),y(t),z(t))$ if $0<|w| <1 $:
\begin{equation}
\begin{gathered}
x(t)=\frac{2c}{w} \sin{\frac{wt}{2}}\cos\Big( \frac{wt}{2}+\alpha \Big),\ \
y(t)=\frac{2c}{w} \sin{\frac{wt}{2}}\sin\Big( \frac{wt}{2}+\alpha \Big), \notag \\
z(t)=wt\cdot \Big\{1+\frac{c^2}{2w^2} \Big[ \Big(1-\frac{\sin(2wt+2\alpha)-\sin{2\alpha}}{2wt}\Big)+ \\
+\Big(1-\frac{\sin(2wt)}{wt}\Big)-\Big(1-\frac{\sin(wt+2\alpha)-\sin{2\alpha}}{2wt}\Big)\Big]\Big\} = \\
=wt\cdot \Big\{1+\frac{c^2}{2w^2} \Big[ \Big(1-\frac{\sin(wt)}{wt}\Big)
+\Big(\frac{1-\cos(2wt)}{wt}\Big) \sin(wt+2\alpha)\Big]\Big\}. \tag{2.8}
\end{gathered}
\end{equation}
In the cases $w=0$ the geodesic curve is the following:
\begin{equation}
x(t)=c\cdot t \cos{\alpha}, \ \ y(t)=c\cdot t \sin{\alpha}, \ \ z(t)=\frac{1}{2} ~ c^2 \cdot t^2 \cos{\alpha} \sin{\alpha}. \tag{2.9}
\end{equation}
The cases $|w|=1$ are trivial: $(x,y)=(0,0), \ z=w \cdot t$.
\begin{definition}
The distance $d^g(P_1,P_2)$ between the points $P_1$ and $P_2$ is defined by the arc length of geodesic curve
from $P_1$ to $P_2$.
\end{definition}
\begin{definition}
 The geodesic sphere of radius $R$ with centre at the point $P_1$ is defined as the set of all points 
 $P_2$ in the space with the condition $d^g(P_1,P_2)=R$. Moreover, we require that the geodesic sphere is a simply connected 
 surface without selfintersection 
 in the $\NIL$ space.
 \begin{rmrk}
 We shall see that this last condition depends on radius $R$.
 \end{rmrk}
 \end{definition}
 \begin{definition}
 The body of the geodesic sphere of centre $P_1$ and of radius $R$ in the $\NIL$ space is called geodesic ball, denoted by $B_{P_1}(R)$,
 i.e., $Q \in B_{P_1}(R)$ iff $0 \leq d^g(P_1,Q) \leq R$.
 \end{definition}
 We proved in\cite{Sz07} the following important theorems:
 \begin{Theorem}[\cite{Sz07}]
 The geodesic sphere and ball of radius $R$ exists in the $\NIL$ space if and only if $R \in [0,2\pi].$
 \end{Theorem}
 \begin{Theorem}[\cite{Sz07}]
 The geodesic $\NIL$ ball $B(S(R))$ is convex in affine-Euclidean sense in our model if and only if $R \in [0,\frac{\pi}{2}]$. 
\end{Theorem}
\subsection{Some properties of geodesic curves and spheres }
In the following, we determine some important properties of geodesic curves and spheres, which we will use in the following sections. 

\begin{enumerate}
\item Consider points $P(x(t),y(t),z(t))$ lying on a sphere $S$ of radius $R$  
centred at the origin. The coordinates of $P$ are given by parameters $(\alpha\in[-\pi, \pi), ~ \theta\in [-\frac{\pi}{2},\frac{\pi}{2}], ~ R>0)$ (see (2.8), (2.9)). 

We obtained directly from the equations (2.8) and (2.9) the following (see \cite{Sz23-1})

\begin{lemma}[\cite{Sz23-1}]
\begin{enumerate}
\item
$$x(t)^2+y(t)^2=\frac{4c^2}{w^2}\sin^2{\frac{wt}{2}},$$
that means, that if $\theta \ne \pm \frac{\pi}{2}$ and $t=R$ is given and 
$\alpha\in[-\pi,\pi)$ then the endpoints $P$ of the geodesic curves lie on a cylinder of radius 
$r=\left| \frac{4c}{w}\sin{\frac{wR}{2}}\right|$ with axis $z$. Therefore, we obtain the following connection between parameters $\theta$ and $R$:
\begin{equation}
r=2\cdot \arcsin \left[  \frac{\sqrt{x^2(R)+y^2(R)}}{2\cdot \cot{\theta}}\right] \frac{1}{\sin{\theta}} \tag{2.10}
\end{equation}
\item
If $\theta = \pm \frac{\pi}{2}$ then the endpoints $P(x(R),y(R),z(R))$ of the geodesics $g^\NIL(x(t),y(t),z(t))$ lie on the 
$z$-axis thus their orthogonal projections onto the $[x,y]$-plane is the origin and 
$x(R)=y(R)=0$, $z(R)=d^g(O,R)=R$.

\item
Moreover, the cross section 
of the spheres $S$ with the plane $[x,z]$ is given by the following system of equation:
     \begin{equation}
     \begin{gathered}
    X(R,\theta)=\frac{2c}{w} \sin{\frac{wR}{2}}=\frac{2\cos{\theta}}{\sin{\theta}} \sin{\frac{R \sin{\theta}}{2}}, \\ 
    Z(R,\theta)=wR+\frac{c^2R}{2w} - \frac{c^2}{2w^2}\sin{wR}= \\ 
    R\sin{\theta}+\frac{R\cos^2{\theta}}{2\sin{\theta}} - \frac{\cos^2{\theta}}{2\sin^2{\theta}}\sin(R\sin{\theta}), \ \ (\theta\in [-\frac{\pi}{2},\frac{\pi}{2}]\setminus\{0\}); \\
    \text{if} \ \theta=0 \ \text{then} \ X(R,0)=R, \ Z(R,0)=0. \tag{2.11}
   \end{gathered}
   \end{equation}
\end{enumerate}
\end{lemma}
\begin{rmrk}
The parametric equations of the geodesic sphere of radius $R$ can be generated from (2.1) by $\NIL$ rotation (see (2.6)).
\end{rmrk}
\item 
We introduce the usual notion of the fibre projection $\mathcal{P}$ that is a projection parallel to fibre lines (parallel to $z$-axis), onto the $[x,y]$ plane.
The image of a point $P$ is the intersection with the $[x,y]$ base plane of the line parallel to fibre line passing through $P$, $\mathcal{P}(P)=P^*$.

Analysed the parametric equations of the geodesic curves $g^\NIL(x(t),y(t),z(t))$ with starting points at the origin we got the following
\begin{lemma}[\cite{Sz23-1}]
If $0 <|w| <1$ for geodesic curve $g^\NIL(x(t),y(t),z(t))$ $(t\in [0,R])$ then 
the fibre projection $\mathcal{P}$ of the geodesic curves (2.8) onto the $[x, y]$ plane is an Euclidean circle arc where it is contained by a circle with equation 
\begin{equation}
\Big(x(t)+\frac{c}{w}\sin{\alpha}\Big)^2+\Big(y(t)-\frac{c}{w}\cos{\alpha}\Big)^2=\Big(\frac{c}{w}\Big)^2=\cot^2{\theta}. \tag{2.12}
\end{equation}
If $w=0$ then fibre projection $\mathcal{P}$ of the geodesic curves $g^\NIL(x(t),y(t),z(t))$ $(t\in [0,R])$ onto the $[x, y]$ plane is a segment with starting point at the origin
where it is contained by the straight line with equation
\begin{equation}
y=\tan{\alpha}\cdot x.\tag{2.13}
\end{equation}
If $w=1$ then the fibre projection $\mathcal{P}$ of the geodetic curves $g^\NIL(x(t),y(t),z(t))$ $(t\in [0,R])$ onto 
the $[x, y]$ plane is the origin.
\end{lemma}
We obtained directly from the equation (2.12) the following
\begin{corollary}[\cite{Sz23-1}]
\begin{enumerate}
\item If we know the equation of the circle that contains the orthogonal projected image $OP^*$ of a geodesic curve segment 
$g^\NIL_{OP}=g^\NIL(x(t),$ $y(t),z(t))$ given in (2.8) $(t\in [0,R])$ onto the $[x,y]$ plane where $0 <|w| <1$ is a known real number (we know the sign of $w$) and the 
coordinates of $P^*=(x(R), y(R),0)$ then the parametric equation of the geodesic curve segment $g^\NIL_{OP}$ is uniquely determined. 
That means that there is {\it one-to-one correspondence between the circle arcs $OP^*$ and the geodesic curve segments $OP$} by the above sense. 
\item If $w=0$ then the fibre projection $\mathcal{P}$ of the geodetic curves is a segment with starting point at the origin
where it is contained by the straight line $y=\tan{\alpha}\cdot x$, therefore in this situation it is 
{\it one-to-one correspondence between the projected image $OP^*$ and the geodesic curve segments $OP$}, too. 
\item If $w=1$ then the fibre projection $\mathcal{P}$ of the geodetic curves 
is the origin so here it is also a {\it one-to-one correspondence between the projected image and the above geodesic curves.}
\end{enumerate}
\end{corollary}
\end{enumerate}
\section{$\NIL$ cylinders and Pappus's hexagon theorem}
\begin{definition}
Let $\cC^i(r)$ be an infinite solid that is ruled by ``side fibre lines" passing 
through the points of an euclidean circle $\cC^b(r)$ of 
radius $r\in \mathbf{R}^+$ ($r$ is the euclidean radius of the cylinder) lying in 
the base plane and centred at the origin.
The images $\cC^i_\bt(r)$ of solid $\cC^i(r)$ by $\NIL$ isometries 
$\bt$ are called {\rm infinite circular cylinders}.
\end{definition}
\begin{rmrk}
In a special case, the Euclidean plane defined by fibre lines fitting to a Euclidean straight line in the base plane is also called a 
``infinite circular cylinder''. In this case, the radius of the ``cylinder" is infinite 
and denoted by $\cC^i_\bt(\infty)$.
\end{rmrk}
The common part of $\cC^i_\bt(r)$ with the base plane is the 
{\it base figure} of $\cC^i_\bt(r)$ that is denoted by $\cC^b_\bt(r)$.
\begin{definition}
A {\rm bounded fibre-like circular cylinder} is an isometric image of a solid
which is bounded by the side surface of an infinite circular cylinder $\cC^i(r)$,
its base figure
$\cC^b(r)$ and its translated copy $\cC^c(r)$ by a fibre translation, 
given by (2.3).
The faces $\cC^b(r)$ and $\cC^c(r)$ are called {\rm cover faces}.
The height (or altitude) of the cylinder is the distance between its cover faces.
\end{definition}
\begin{figure}[ht]
\centering
\includegraphics[width=10cm]{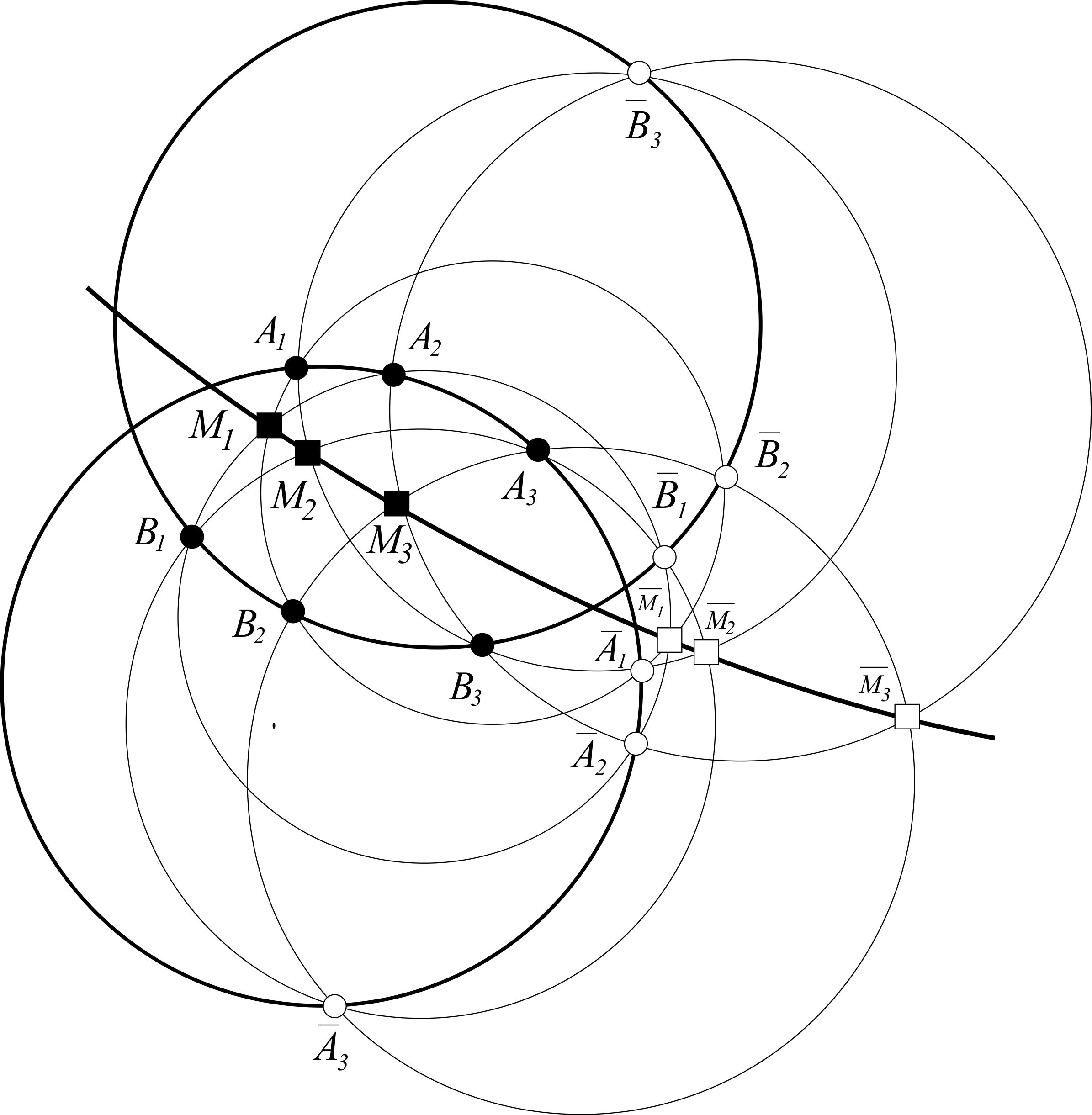}
\caption{Circle-like Pappus's hexagon theorem}
\label{}
\end{figure}
The straightforward consequence of the formulas (2.2-3), (2.6), (2.8-9) and 
the Definitions 3.1-2 the following 
\begin{lemma}
If the euclidean radius of the circle $\cC^b(r)$ centred at the origin lying in the base 
plane is $r$ then the corresponding infinite fibre-like cylinder $\cC^i(r)$
is a circual cylinder in Euclidean and $\NIL$ sense, too (see Fig.~2-3). \qquad \qquad $\square$
\end{lemma}
Let us denote the image of the cylinder $\cC^i(r)$ at the $\NIL$ translation 
$\bT$ (see (2.3))
by $\cC^i_\bT(r)$ where the translation is given by parameters 
$x^0=1,x^1,x^2,x^3$ and denote
its common part with the base plane by $\cC^b_\bT(r)$.

Using Lemma 2.9, Lemma 2.11 and formulas (2.2), (2.3)  we obtain the following
\begin{lemma}
$\cC^b_\bT(r)$ is an euclidean circle of radius $r$.  \qquad \qquad $\square$
\end{lemma}
\subsection{Pappus-like theorem for circles in inversive plane}

The classical M\"obius plane is the Euclidean plane supplemented by a single point at infinity. It is also called the inversive plane because it is closed under 
inversion with respect to any circle. The following theorem can be proven using the usual inversion and elementary geometric theorems and relations. 
\begin{Theorem}
Consider six arbitrary points on the inversive plane 
$A_1, A_2, A_3,$ $B_1, B_2, B_3$. We consider two circles (can also be a straight line, i.e. a circle with infinite radius) 
$K_{A_1A_2A_3}$ passing through the points $A_1,A_2,A_3$ and $K_{B_1B_2B_3}$ containing $B_1, B_2, B_3$. Let a point $\overline{A}_1$ be on the circle $K_{A_1A_2A_3}$. 
Let $\overline{B}_2$ be the point where the circle $K_{A_1B_2\overline{A}_1}$ intesects the circle $K_{B_1B_2B_3}$. Moreover, draw the circle $K_{B_2A_3\overline{B}_2}$ and let 
$\overline{A}_3=K_{B_2A_3\overline{B}_2} \cap K_{A_1A_2A_3}$ where $\overline{A}_3 \ne \overline{B_2}$. Let $\overline{B}_1$ be the from $B_1$ different 
point where the circle $K_{A_3B_1\overline{A}_3}$
intesects the circle $K_{B_1B_2B_3}$. Then, we consider the circle $K_{B_1A_2\overline{B}_1}$ and let $\overline{A}_2$ be the other point where it intesects the circle $K_{A_1A_2A_3}$. 
Let $\overline{B}_3$ be the intersection point of the circles $K_{A_2B_3\overline{A}_2}$ and $K_{B_1B_2B_3}$ which is different from $B_3$. 
\begin{enumerate}
\item Thus, $\overline{A}_1 \in K_{B_3A_1\overline{B}_3}$. 
\item Let $M_1$ and $\overline{M}_1$ be the points where the circles $K_{A_1B_2 \overline{A}_1 \overline{B}_2}$ and $K_{A_2B_1 \overline{A}_2 \overline{B}_1}$ meet, 
and $M_2$ and $\overline{M}_2$ where 
$K_{A_1B_3\overline{A}_1,\overline{B}_3}$ and $K_{A_3B_1\overline{A}_3,\overline{B}_1}$ meet, moreover, $M_3$ and $\overline{M}_3$ 
where $K_{A_3B_2\overline{A}_3,\overline{B}_2}$ and $K_{A_2B_3\overline{A}_2,\overline{B}_3}$ meet. 
Then $M_1, M_2, M_3$, $\overline{M}_1$, $\overline{M}_2$ and $\overline{M}_3$ are concyclic. 
\end{enumerate} \qquad $\square$
\end{Theorem}
\begin{rmrk}
If the points $\overline{A}_i$ and $\overline{B}_i$ $i \in \{1,2,3 \}$ coincide, then using an
inversion we directly obtain the classical Pappus hexagon configuration.
\end{rmrk}
\subsection{Pappus-like theorem for $\NIL$ cylinders}
Let us formulate the next  theorem similar to Theorem 3.6 for $\NIL$ infinite cylinders.
\begin{Theorem}
Consider six arbitrary points on the base plane of the considered model of the $\NIL$ geometry. 
$A_1, A_2, A_3,$ $B_1, B_2, B_3$. We consider two circular cylinders (can also be a euclidean plane, i.e. a cylinder with infinite radius) 
$\cC^i_{A_1A_2A_3}$ passing through the points $A_1,A_2,A_3$ and $\cC^i_{B_1B_2B_3}$ containing $B_1, B_2, B_3$. Let a fibre line $\overline{a}^f_1$ be on the cylinder 
$\cC^i_{A_1A_2A_3}$ and let its common point be with the base plane $\overline{A}_1$ 
Let $\overline{b}^f_2$ be the fibre line where the cylinder $\cC^i_{A_1B_2\overline{A}_1}$ intesects the cylinder $\cC^i_{B_1B_2B_3}$ and we denote 
its intersection point with the base plane by $\overline{B}_2$. Moreover, consider the cylinder $\cC_{B_2A_3\overline{B}_2}$ and let 
$\overline{a}^f_3=\cC^i_{B_2A_3\overline{B}_2} \cap \cC^i_{A_1A_2A_3}$ where $\overline{A}_3 \ne \overline{B}_2$. Let $\overline{b}^f_1$ be a fibre line 
where the cylinder $\cC^i_{A_3B_1\overline{A}_3}$
intesects the cylinder $\cC^i_{B_1B_2B_3}$ moreover its common point $\overline{B}_1$ with the base plane is different from $B_1$. 
Then, we consider the cylinder $\cC^i_{B_1A_2\overline{B}_1}$ and let $\overline{a}^f_2$ be the other fibre line where it intesects the cylinder $\cC^i_{A_1A_2A_3}$ and 
let $\overline{A}_2$ the intersection pont of $\overline{a}^f_2$ with the base plane of the model. 
Let $\overline{b}^f_3$ be the intersection fibre line of the cylinders $\cC^i_{A_2B_3\overline{A}_2}$ and $\cC^i_{B_1B_2B_3}$ which is different from the fibre line through $B_3$, 
moreover, we denote its intersection point with the base plane by $\overline{B}_3$.
\begin{figure}[ht]
\centering
\includegraphics[width=12cm]{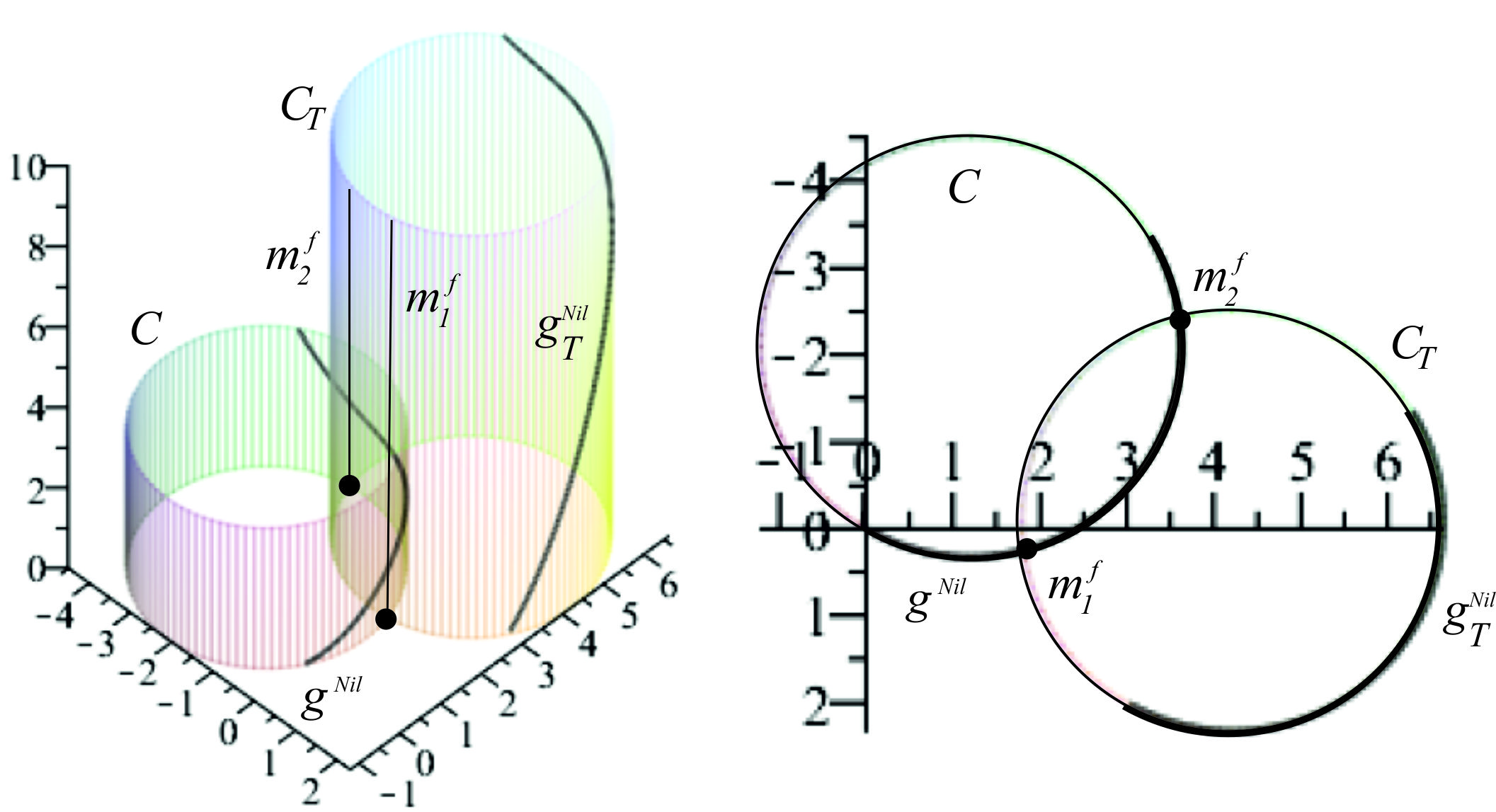}
\caption{Two congruent helix-like geodesic curves and the congruent cylinders containing them.}
\label{}
\end{figure}
\begin{enumerate}
\item Thus, $\overline{a}^f_1 \in \cC_{B_3A_1\overline{B}_3}$. 
\item Let $m^f_1$ and $\overline{m}^f_1$ be the fibre lines where the cylinders $\cC^i_{A_1B_2 \overline{A}_1 \overline{B}_2}$ and $\cC^i_{A_2B_1 \overline{A}_2 \overline{B}_1}$ meet, 
$m^f_2$ and $\overline{m}^f_2$ where 
$\cC^i_{A_1B_3\overline{A}_1,\overline{B}_3}$ and $\cC^i_{A_3B_1\overline{A}_3,\overline{B}_1}$ meet, moreover, $m_3^f$ and $\overline{m}_3^f$ 
where $\cC^i_{A_3B_2\overline{A}_3,\overline{B}_2}$ and $\cC^i_{A_2B_3\overline{A}_2,\overline{B}_3}$ meet. 
Then the fibre lines $m_1^f, m_2^f, m_3^f$, $\overline{m}_1^f, \overline{m}_2^f$ and $\overline{m}_3^f$ lie in a cylinder. 
\end{enumerate} \qquad $\square$
\end{Theorem}
Based on the Lemmas 2.11, 3.4, 3.5 and Corollary 2.12 the infinite circular cylinders determine the equivalence classes of geodesic curves up to $\NIL$ congruences, 
i.e., the Euclidean radii of the 
$\NIL$ infinite circular cylinders characterizes the equivalence classes of geodesic curves. 
Based on these, using the notations of Theorem 3.7 we can state the following theorems:
\begin{Theorem}
There exists uniquely equivalence class of $\NIL$ geodesics that intersects the fibers $m_1^f, m_2^f, m_3^f, \overline{m}_1^f, \overline{m}_2^f$ 
and $\overline{m}_3^f$ \qquad \qquad $\square$
\end{Theorem}
It would not be easy to visualize the many cylinders that occur in the theorem, 
so we only show the intersection of two geodesics belonging to different equivalence classes in Fig.~3.
\begin{figure}[ht]
\centering
\includegraphics[width=12cm]{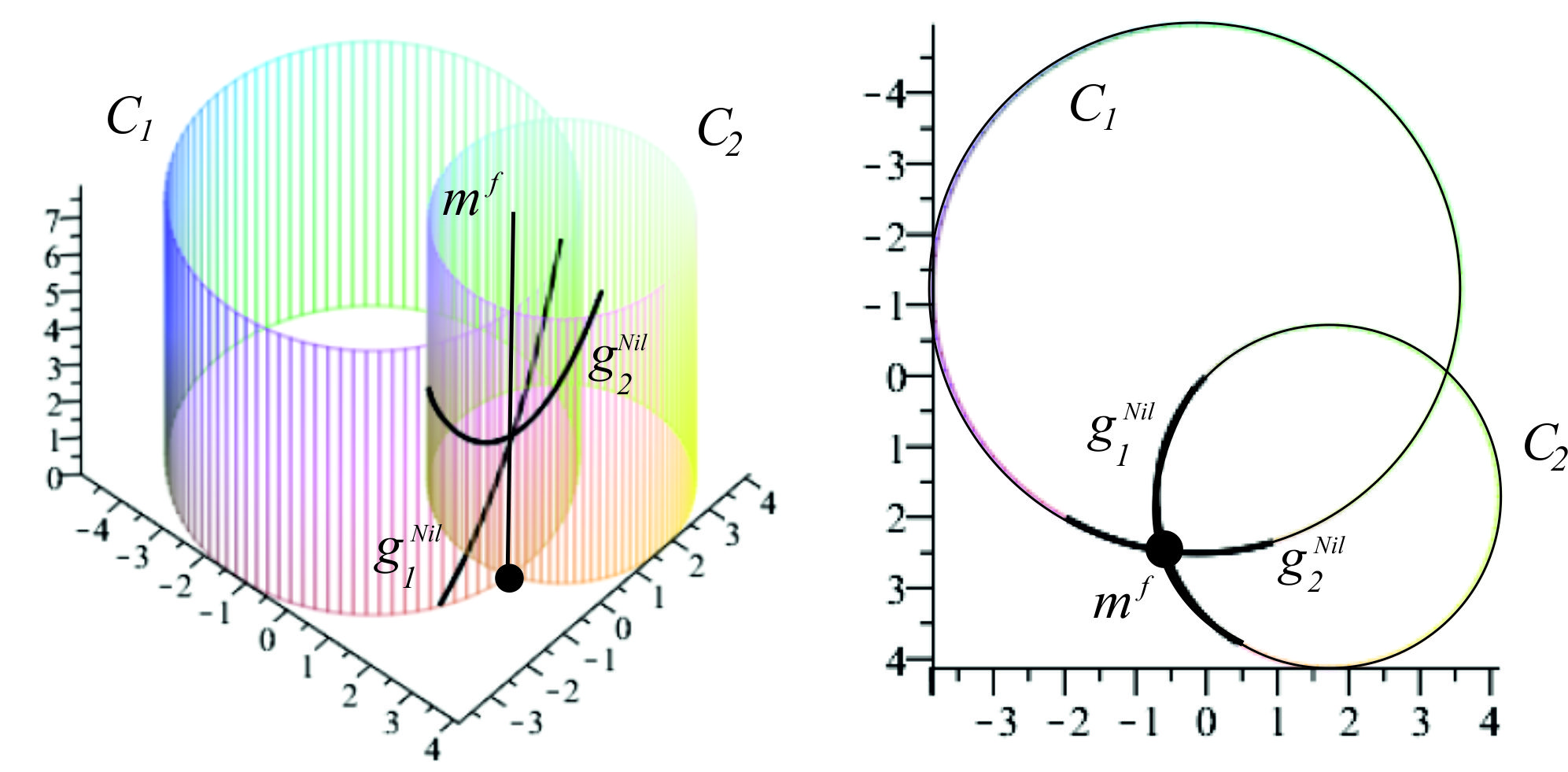}
\caption{Two geodesics belonging to different equivalence classes intersect each other on a fiber.}
\label{}
\end{figure}
\begin{Theorem}
Let $g^\NIL_{m_1^fm_2^fm_3^f}\subset \cC_{m_1^fm_2^fm_3^f}$ be an given geodesic curve (see Theorem 3.7). 
There always exist representation elements of the equivalence classes of $\NIL$ geodesics related to the Theorem 3.7
so that their intersection points 
\begin{enumerate}
\item $M_1 = g^\NIL_{A_1B_2 \overline{A}_1 \overline{B}_2} \cap g^\NIL_{A_2B_1 \overline{A}_2 \overline{B}_1} \in m_1^f$,
\item $M_2 = g^\NIL_{A_1B_3\overline{A}_1,\overline{B}_3} \cap g^\NIL_{A_3B_1\overline{A}_3,\overline{B}_1} \in m_2^f$,
\item $M_3 = g^\NIL_{A_3B_2\overline{A}_3,\overline{B}_2} \cap g^\NIL_{A_2B_3\overline{A}_2,\overline{B}_3} \in m_3^f$
\end{enumerate}
exist and lie on $g^\NIL_{m_1^fm_2^fm_3^f}\subset \cC_{m_1^fm_2^fm_3^f}$. \qquad \qquad $\square$
\end{Theorem}
As it is known, the Pappus's hexagon theorem plays a key role in the axiomatic structure of projective geometry and
now we proved that an analogous theorem is also true in $\NIL$ geometry. Therefore, it is interesting to 
examine the role of this theorem in the structure of $\NIL$ geometry and in general 
in other non-constant curvature Thurston geometries.
Detailed studies are the objective of ongoing research.
\medbreak


\begin{thebibliography}{12}
%
\bibitem{CaMoSpSz}
{Cavichioli,~A.~--~Moln\'ar,~E.~--~Spaggiari,~F.~--~Szirmai,~J.,}
Some tetrahedron manifolds with $\SOL$ geometry.
\textit{J. Geom.,} {\bf 105/3}, 601-614 (2014).
%
\bibitem{CsSz16}
{Csima,~G.~--~Szirmai,~J.},
{Interior angle sum of translation and geodesic triangles in $\SLR$ space.}
\emph{Filomat}, {\bf 32/14}, (2018) 5023--5036.
%
\bibitem{Cs-Sz23}
{Csima,~G.~--~Szirmai,~J.,}
Translation-like isoptic surfaces and angle sums of translation triangles in $\NIL$ geometry.
\emph{Results Math.}, (2023) 78:194, DOI: 10.1007/s00025-023-01961-z.
%
\bibitem{Cs-Sz25}
{Csima,~G.~--~Szirmai,~J.,}
Translation-like Apollonius and triangular surfaces in non-constant curvature Thurston geometries.
\emph{Results Math.}, {\bf 80}, article number 190, (2025), DOI: 10.1007/s00025-025-02503-5.
%
\bibitem{M97}
{Moln{\'a}r,~E.,}
The projective interpretation of the eight 3-di\-men\-sional homogeneous geometries.
\emph{Beitr. Algebra Geom.,}
{\bf38} No.~2, 261--288, (1997).
%
\bibitem{M06}
{Moln{\'a}r,~E.:}
On projective models of Thurston geometries, some relevant notes
on $\NIL$ orbifolds and manifolds. \emph{Sib. Electron.  Math. Izv.},
{\bf 7} (2010), 491--498,
http://mi.mathnet.ru/semr267
%
\bibitem{MSz06}
{Moln{\'a}r,~E.~--~Szirmai,~J.,}
{On $\NIL$ crystallography,}
{\it Symmetry Cult. Sci.,}
{\bf 17/1-2} (2006), 55--74.
%
\bibitem{MSzV}
Moln{\'a}r,~E.~--~Szirmai,~J.~--~Vesnin,~A.,
Projective metric realizations of cone-manifolds with singularities along 2-bridge knots and links.
{\it J. Geom.,}  {\bf 95}, 91-133 (2009).
%
\bibitem{MSz}
{Moln{\'a}r,~E.~--~Szirmai,~J.,}
Symmetries in the 8 homogeneous 3-geometries.
\textit{Symmetry Cult. Sci.,}
{\bf 21/1-3}, 87-117 (2010).
%
\bibitem{MSz12}
{Moln{\'a}r,~E.~--~Szirmai,~J.,}
Classification of $\SOL$ lattices.
\textit{Geom. Dedicata,}
{\bf 161/1}, 251-275 (2012).
%
\bibitem{PSSz10}
{Pallagi,~J.~--~Schultz~B.~--~Szirmai,~J.},
{Equidistant surfaces in $\NIL$ space,}
\emph{Stud. Univ. Zilina, Math. Ser.,}
{\bf 25}, 31--40 (2011).
%
\bibitem{PS14}
{Papadopoulos,~A.~--~Su,~W.},
On hyperbolic analogues of some classical theorems in spherical geometry.
(2014), hal-01064449.
%
\bibitem{S}
Scott,~P.,
The geometries of 3-manifolds. {\it Bull. London Math. Soc.}  {\bf 15}, 401--487 (1983).
%
\bibitem{Sz07}
Szirmai,~J.,
The densest geodesic ball packing by a type of $\NIL$ lattices.
{\it Beitr. Algebra Geom.} {\bf 48}(2) (2007), 383--398.
%
\bibitem{Sz14-1}
Szirmai,~J.,
A candidate to the densest packing with equal balls in the Thurston geometries. 
{\it Beitr. Algebra Geom.,} {\bf 55}(2) (2014), 441--452.
%
\bibitem{Sz14-2}
{Szirmai,~J.},
Simply transitive geodesic ball packings to $\mathbf{S^2\times R}$ space groups generated by glide reflections,
{\emph {Ann. Mat. Pur. Appl.}}, {\bf 193/4} (2014), 1201-1211, DOI: 10.1007/s10231-013-0324-z.
%
\bibitem{Sz18}
Szirmai,~J.,
{$\NIL$ geodesic triangles and their interior angle sums},  
{\it Bull. Braz. Math. Soc. (N.S.),} {\bf 49} (2018) 761--773, DOI: 10.1007/s00574-018-0077-9.
%
\bibitem{Sz19}
Szirmai,~J.,
Bisector surfaces and circumscribed spheres of tetrahedra derived
by translation curves in $\SOL$ geometry. 
{\it New York J. Math.,} {\bf 25}, 107--122 (2019).
%
\bibitem{Sz22-2}
Szirmai,~J.,
Apollonius surfaces, circumscribed spheres of tetrahedra, Menelaus' and Ceva's theorems in $\SXR$ and $\HXR$ geometries. 
{\it Q. J. Math.,} {\bf 73}, (2022), 477-494, DOI: 10.1093/qmath/haab038.
%
\bibitem{Sz23-1}
Szirmai,~J.,
On Menelaus' and Ceva's theorem in $\NIL$ geometry. 
{\it Acta Univ. Sapientiae Math.,} {\bf 15/ 1}, (2023) 123--141, DOI: 10.2478/ausm-2023-0008.
%
\bibitem{Sz23-2}
{Szirmai,~J.,}
Classical Notions and Problems in Thurston Geometries,
\emph{International Electronic Journal of Geometry},
{\bf 16} No.2 (2023), 608--643, doi: 10.36890/IEJG.1221802.
%
\bibitem{Sz24}
{Szirmai,~J.},
Fibre-like cylinders, their packings and coverings in $\SLR$ space,
\emph{Results Math.}, {\bf{79}} article number 123 (2024), doi: 10.1007/s00025-024-02152-0.
%
\bibitem{Sz26}
{Szirmai,~J.},
Menelaus' and Ceva's theorems for translation triangles in Thurston geometries, 
\emph{Results Math.}, (2026) (to appear), arXiv: 2506.01354.
%
\bibitem{Sz26-1}
{Szirmai,~J.},
Desargues's and Pappus's hexagon theorems on translation-type surfaces in Thurston geometries, 
\emph{Submitted manuscript}, (2026), arXiv: 2603.00019.
%
\bibitem{T}
Thurston,~W.~P. (and Levy,~S. editor),
{\it Three-Dimensional Geometry and Topology}.  Princeton University Press,  Princeton, New Jersey, vol. {\bf 1} (1997).
%
\bibitem{VSz19}
{Vr\'anics,~A.~--~Szirmai,~J.},
Lattice coverings by congruent translation
balls using translation-like bisector
surfaces in Nil Geometry.
\emph{KoG}, {\bf 23}, 6-17 (2019).
%
\bibitem{YSz24-2}
{Yahya, A.~--~Szirmai, J.,}
Geodesic ball packings generated by rotations and monotonicity behavior of their densities in $\mathbf{H}^2\!\times\!\mathbf{R}$ space.
\emph{Results Math.}, {\bf{80}}, 123 (2025), doi: 10.1007/s00025-025-02430-5.
%
\end{thebibliography}
\end{document}